%fcp.tex: 
%%a Plain TeX file by  Doron Zeilberger and Noam Zeilberger (x pages)

%begin macros

\baselineskip=14pt
\parskip=10pt
\def\halmos{\hbox{\vrule height0.15cm width0.01cm\vbox{\hrule height
  0.01cm width0.2cm \vskip0.15cm \hrule height 0.01cm width0.2cm}\vrule
  height0.15cm width 0.01cm}}

\def\1{{\overline{1}}}
\def\2{{\overline{2}}}
\parindent=0pt
\overfullrule=0in

\def\frac#1#2{{#1 \over #2}}
%\headline={\rm  \ifodd\pageno  \RightHead  \else  \LeftHead  \fi}
%\def\RightHead{\centerline{
%Title
%}}
%\def\LeftHead{ \centerline{Doron Zeilberger}}
%end macros
\centerline
{\bf Two Questions about the Fractional Counting of Partitions }
\bigskip
\centerline
{\it Doron ZEILBERGER and Noam ZEILBERGER}
\bigskip

{\bf Abstract}: We recall the notion of {\it fractional enumeration} and immediately focus on
the fractional counting of integer partitions, where each partition gets `credit' equal to the
reciprocal of the product of its parts. We raise two intriguing questions regarding this count,
and for each such question we are pledging a \$100 donation to the OEIS, in honor of the first solver.
In this revised version we announce that both questions have been answered.
The first was first answered by Will Sawin and 
the second was answered by Christopher Ryba. A donation to the OEIS, of \$100 each, in their honor
has been made.

{\bf Update Oct. 31, 2018}: Will Sawin has solved the first question. A donation to the OEIS, in
his honor, has been made. His proof is given as an appendix to this revised version.]

{\bf Update 10am, Nov. 2, 2018}: Fernando Chamizo has independently solved the first question, with an
even shorter rendition. We give a reference to it in this revised version.]

{\bf Update 2:30pm, Nov. 2, 2018}:  Gjergji Zaimi kindly pointed out that Question 1 was posed and answered
completely by D.H. Lehmer in Theorem 3 of his article 
``On reciprocally weighted partitions", Acta Arithmetica XXI(1972), p. 379-388, that is available here: \hfill\break
 {\tt http://matwbn.icm.edu.pl/ksiazki/aa/aa21/aa21123.pdf} \quad .

{\bf Update  Nov. 8, 2018}:  Christopher Ryba solved  Question 2 . His elegant article
appeared here:\hfill\break
{\tt  https://arxiv.org/abs/1811.03440} \quad .

{\bf Preface: There are many ways to enumerate}

Naive counting of combinatorial sets
counts by $1$. Generatingfunctionology counts by $z^{stat}$, where $z$ is an {\it indeterminate} (i.e. symbolic), and
$stat$ is a certain statistic of interest. Statistical mechanics also `counts' by $z^{stat}$, but now $z$ is
a continuous real (or complex) variable of physical significance. Sieve theory `counts' by $\pm 1$.

The general scenario of naive enumeration is a sequence of combinatorial sets, $A_n$, naturally indexed by a non-negative integer $n$,
and one wants a formula, or at least an efficient algorithm, to compute the number of elements of $A_n$, denoted by
$|A_n|$.

For example, 

$\bullet$ The sequence `set of subsets of  $\{1, \dots, n \}$', where $|A_n|=2^n$. \hfill\break
It is:  {\tt  https://oeis.org/A000079} \quad .

$\bullet$ The sequence `set of permutations of  $\{1, \dots, n \}$', where $|A_n|=n!$. \hfill\break
 It is:  {\tt https://oeis.org/A000142} \quad .

$\bullet$ The sequence of integer partitions of  $n$, where $|A_n|=p(n)$. \hfill\break
It is: { https://oeis.org/A000041} \quad .

There are no `nice' formulas for $p(n)$, but there exist efficient algorithms.
One not so efficient `formula'  is 

\centerline{`the coefficient of $q^n$ in $1/((1-q)(1-q^2) \cdots (1-q^n))$' \quad .}

The general scenario of Generatingfunctionology
enumeration is a sequence of combinatorial sets, $A_n$, naturally indexed by a non-negative integer $n$,
and a certain `statistic' defined on its objects $s \rightarrow stat(s)$,
and one wants a formula, or at least an efficient algorithm, to compute the sequence
of {\it weight-enumerators}, $|A_n|_z:=\sum_{s \in A_n} z^{stat(s)}$.

For example, 

$\bullet$ The sequence  `set of subsets of  $\{1, \dots, n \}$', where the statistic is `cardinality', and we have
$|A_n|_z=(1+z)^n$.

$\bullet$ The sequence `set of permutations of  $\{1, \dots, n \}$', where the statistic is `number of inversions', and we have
$|A_n|_z=1 \cdot (1+z) \cdot (1+z+z^2) \cdots (1+z+ \dots +z^{n-1})$ \quad .

$\bullet$ The sequence of integer partitions of  $n$, where the statistic is `largest part',
for which, once again, there is no `nice' formula, but there exist efficient algorithms.
One, not so efficient `formula'  is the coefficient of $q^n$ in $1/((1-qz)(1-q^2z) \cdots (1-q^n z))$

{\bf Signed Counting}

Let $\mu(n)$ be $1$ if $n$ is a product of an even number of distinct primes, $-1$ if it is a product of 
an odd number of distinct primes, and $0$ otherwise. The signed enumeration of the set $A(n):=\{1, \dots , n\}$,
whose naive count is $n$, is 
$$
M(n):=\sum_{i=1}^{n} \mu(i) \quad .
$$

{\bf Exercise}: Find a closed form formula for $M(n)$ that would entail its asymptotic behavior. If you can't,
prove at least that it is $O(n^{\frac{1}{2}+\epsilon})$ for any $\epsilon>0$. Failing this, prove at least that it is
$O(n^{0.99999999999999999999})$.

{\bf Functional Counting} 

The function $f(x):=z^x$, that gives the generatingfunctionology count could be replaced by {\it any} function. So let's define
$$
a_f(n):=\sum_{s \in A_n} \, f(stat(s)) \quad .
$$
If $f$ is a power $f(x):=x^k$, one gets the numerator of the $k$-th moment of $stat$.

{\bf Fractional Counting} 

This corresponds to the case $f(x) := \frac{1}{x}$.
Fractional counting of lattice paths is dealt with in [FGZ].
Another natural form of fractional counting was described by Baez and Dolan [BD], called `homotopy cardinality' (or `groupoid cardinality').
The homotopy cardinality
$$
|X| = \sum_{[x] \in X/\sim}\frac{1}{|{\rm Aut}(x)|}
$$
of a groupoid $X$ counts the number of isomorphism classes of objects in $X$, but where each class is weighted inversely by the size of the automorphism group of any representative element.

{\bf Maple package} 

This article is accompanied by a Maple package, {\tt FCP.txt}, obtainable from the front of this article:

{\tt http://sites.math.rutgers.edu/\~{}zeilberg/mamarim/mamarimhtml/fcp.html} \quad .

That page also contains sample input and output files, as well as a nice picture.

{\bf Fractional Counting of Partitions}

From now we will focus on fractional counting of
integer partitions where each partition gets `credit' the reciprocal of the product of its parts.

{\bf Definition:} Let $b(n)$ be defined by
$$
b(n):= \sum_{
{{p_1+ \dots + p_k=n}
\atop
 {p_1 \geq p_2 \geq \dots \geq p_k >0}
}
}
\frac{1}{p_1 p_2 \dots p_k} \quad .
$$

Let's spell out the first few terms of the sequence of fractions $b(n)$
$$
b(1)\, = \, \frac{1}{1} =1 \quad ,
$$
$$
b(2) \, = \, \frac{1}{1 \cdot 1} + \frac{1}{2} = \frac{3}{2} \quad ,
$$
$$
b(3) \,  = \, \frac{1}{1 \cdot 1 \cdot 1} + \frac{1}{2\cdot 1} + \frac{1}{3} = \frac{11}{6} \quad .
$$

For the first $100$ terms see

{\tt http://sites.math.rutgers.edu/\~{}zeilberg/tokhniot/oFCP1.txt} \quad .

{\bf How to compute $b(n)$ for many $n$?}

Let's recall one of the many ways to compute a table of $p(n)$, the naive count of the
set of partitions of $n$, for many values of $n$.
It is not the most efficient way, but the one easiest to adapt to the computation of the {\it fractional} count, that
we called $b(n)$.

Let $p(n,k)$ be the number of partitions of $n$ whose largest part is $k$. Once we know $p(n,k)$ for
$1 \leq k \leq n \leq N$ we would, of course, know $p(n)$ for $1 \leq n \leq N$, since
$$
p(n) \, = \, \sum_{k=1}^{n} p(n,k) \quad .
$$

$p(n,k)$ may be computed, recursively, via the {\it dynamical programming} recurrence
$$
p(n,k)= \sum_{k'=1}^{k} p(n-k,k') \quad ,
$$
since removing the largest part, $k$, of a partition of $n$ results in a partition of $n-k$ whose largest part is $\leq k$.

A more compact recursion, without the $\sum$, is obtained as follows.
Replace, in the above equation,  $n$ and $k$ by $n-1$ and $k-1$ respectively, to get
$$
p(n-1,k-1)= \sum_{k'=1}^{k-1} p(n-k,k') \quad .
$$
Subtracting gives
$$
p(n,k) \, - \, p(n-1,k-1) \, = \, p(n-k,k) \quad ,
$$
leading to the recurrence
$$
p(n,k) \, = \,  p(n-1,k-1) \, + \, p(n-k,k) \quad .
$$

We now proceed analogously.

Let $b(n,k)$ be the fractional count of the set of partitions of $n$ whose largest part is $k$. Once we know $b(n,k)$ for
$1 \leq k \leq n \leq N$ we would, of course, know $b(n)$ for $1 \leq n \leq N$, since
$$
b(n) \, = \, \sum_{k=1}^{n} b(n,k) \quad .
$$

$b(n,k)$ may be computed via the {\it dynamical programming} recurrence
$$
b(n,k)= \frac{1}{k} \sum_{k'=1}^{k} b(n-k,k') \quad ,
$$
since removing the largest part, $k$, of a partition of $n$ results in a partition of $n-k$ whose largest part is $\leq k$.

Multiplying both sides by $k$ yields
$$
k\, b(n,k) \, = \, \sum_{k'=1}^{k} b(n-k,k') \quad .
$$

Replacing $n$ and $k$ by $n-1$ and $k-1$ respectively, yields
$$
(k-1) \, b(n-1,k-1)= \sum_{k'=1}^{k-1} b(n-k,k') \quad .
$$
Subtracting gives
$$
k \, b(n,k) \, - \, (k-1)\, b(n-1,k-1) \, = \, b(n-k,k) \quad ,
$$
leading to the recurrence
$$
b(n,k) \, = \,  \frac{k-1}{k} \, b(n-1,k-1) \, + \, \frac{1}{k} \, b(n-k,k) \quad ,
$$
with  the boundary conditions $b(n,1)=1$ and $b(n,k)=0$ if $k>n$.

Procedure {\tt bnk(n,k)} implements this recurrence in the Maple package {\tt FCP.txt}, and 
{\tt bnkF(n,k)} is the much faster floating-point version.

We believe that the following fact is easy to prove.

{\bf Fact:} $C:= \lim_{n \rightarrow \infty } \frac{b(n)}{n} \quad$ exists.

The convergence is rather slow. Here are the values of $\frac{b(n)}{n}$ for $15000-10\leq n \leq 15000$.
$$
0.5611411658 \, ,  \, 0.5611411846 \, ,  \, 0.5611412033 \, , \, 0.5611412220 \, ,  \, 0.5611412407 \, , \,  0.5611412594, 
$$
$$
0.5611412781 \, , \, 0.5611412968 \, , \, 0.5611413156 \, ,  \, 0.5611413344  \, , \,  0.5611413530 \quad .
$$

One of us (DZ) is pledging a \$100 donation to the OEIS in honor of the first person to answer  the following question.

{\bf Question 1}: Identify $C$ in terms of known mathematical constants. In particular, is $C=e^{-\gamma}$?
Here $\gamma$ is Euler's constant. Note that $e^{-\gamma}=0.5614594835668851698\dots$ \quad .

{\bf Added Oct. 31, 2018}: We are almost sure that $C$ is indeed $e^{-\gamma}$. If you use the `ansatz'
$$
\frac{b(n)}{n} \, = \, c_0 + \frac{c_1}{n} \quad ,
$$
and plug-in $n=14999$ and $n=15000$, you would get two equations in the two unknowns $c_0$ and $c_1$, whose solution is
$$
   \{c_0 =0.5614203344, c_1 = -4.184721000\} \quad,
$$
agreeing to four decimal places.

[Update Oct. 31, 2018: Will Sawin has solved  Question 1. A donation to the OEIS, in
his honor, has been made. His proof is given as an appendix to this article.]

[Update 10:30am, Nov. 2, 2018: Fernando Chamizo has independently solved the first question, with an
even shorter rendition. It is available, with kind permission of Professor Chamizo, from \hfill\break
{\tt http://sites.math.rutgers.edu/\~{}zeilberg/mamarim/mamarimPDF/fcpFernandoChamizo.pdf } \quad .]

[Update 2:30pm, Nov. 2, 2018:  Gjergji Zaimi kindly pointed out that Question 1 was posed and answered
completely by D.H. Lehmer in Theorem 3 of his article  
``On reciprocally weighted partitions", Acta Arithmetica XXI(1972), p. 379-388, \hfill\break
{\tt http://matwbn.icm.edu.pl/ksiazki/aa/aa21/aa21123.pdf } ]

We also noticed, numerically, that for each real $0 <x<1$
$$
f(x) := lim_{n \rightarrow \infty} \,\, b(n,\lfloor nx \rfloor) \quad,
$$
exists, and defines a nice, decreasing function. To see an approximation, using $n=2000$, see

{\tt http://sites.math.rutgers.edu/\~{}zeilberg/tokhniot/picsFCP/fcp1.html} \quad \quad .

We are also  pledging a \$100 donation to the OEIS in honor of the first person to answer the following question.

{\bf Question 2}: Find a differential equation satisfied by $f(x)$, and if possible, an explicit expression in terms of
known functions. 

Note that
$$
C \, = \, \int_{0}^{1} \, f(x) \, dx \quad ,
$$
so an answer to Question 2 may settle Question 1.

{\bf Update  Nov. 8, 2018}:  Christopher Ryba solved  Question 2 . His elegant article
appeared here:\hfill\break
{\tt  https://arxiv.org/abs/1811.03440} \quad .

{\bf A  natural approach}

The asymptotic expression for the naive counting of partitions, the famous partition function $p(n)$, { https://oeis.org/A000041},
is the subject of the celebrated Hardy-Ramanujan-Rademacher formula (See [A], Ch. 6). The proof uses the {\it Circle method}, that uses
the Euler generating function
$$
\sum_{n=0}^{\infty} p(n)\, q^n \, = \, \frac{1}{(1-q)(1-q^2)(1-q^3) \cdots } \quad .
$$
This enables one to express $p(n)$ as a contour integral.

The analogous generating function for our sequence of interest is obviously

$$
\sum_{n=0}^{\infty} b(n)\, q^n \, = \, \frac{1}{(1-q)(1-\frac{q^2}{2})(1-\frac{q^3}{3}) \cdots } \quad .
$$
It is possible that a similar proof (possibly easier, since we do not want the full asymptotics only the leading term)
would solve Question 1.

Regarding $b(n,k)$, we obviously have
$$
\sum_{n=0}^{\infty} b(n,k)\, q^n \, = \, \frac{q^k/k}{(1 \, - \, q)\, (1 \, - \, \frac{q^2}{2})\, (1 \, - \, \frac{q^3}{3}) \cdots (1 \,- \,\frac{q^k}{k})} \quad .
$$
Once again this implies a certain contour integral expression for $b(n,k)$ that may lead to an answer to Question 2.

{\bf A much easier question}

If you look at partitions in {\bf frequency notation} $1^{a_1} \dots n^{a_n}$ and give each of them `credit'
$$
\frac{1}{ (1^{a_1} a_1!) \cdots (n^{a_n} a_n!)} \quad ,
$$ 
and define $c(n)$ as the sum of these over all partitions of $n$, then
we have
$$
\sum_{n=0}^{\infty} c(n) q^n = e^{q/1} \, e^{q^2/2} \, e^{q^3/3} \cdots = exp\, \left ( \sum_{i=1}^{\infty} \frac{q^i}{i} \right )
\, = \,
exp(-log(1-q))= \frac{1}{1-q} \quad ,
$$
and it follows that $c(n)=1$ for all $n$.

Recall that $\frac{n!}{ (1^{a_1} a_1!) \cdots (n^{a_n} a_n!)}$ is famously the number of permutations of $n$ whose cycle structure is
$1^{a_1} \dots n^{a_n}$, i.e. the number of permutations of $n$ whose expression into a product of disjoint cycles has $a_1$ fixed points,
$a_2$ cycles of length $2$ etc. Recall the easy proof of that fact. If you order each of the cycle-lengths, then the number of
ways of distributing $n$ `balls'  into the $a_1+ \dots +a_n$ `boxes' is the multinomial coefficient 
$\frac{n!}{ 1!^{a_1} 2!^{a_2} \cdots n!^{a_n}}$.
Since the order of cycles of the same length does not matter, we have to divide by $a_1! \cdots a_n!$, getting
$\frac{n!}{(a_1! 1!^{a_1}) (a_2! 2!^{a_2}) \cdots (a_n! n!^{a_n})}$. But each box with, say, $r$ `balls', can be arranged
into $(r-1)!$ cycles hence the desired number has to be multiplied by $1!^{a_2} 2!^{a_3} \cdots (n-1)!^{a_n}$ proving this fact,
that goes back, at least, to Cauchy and Cayley.  Since the total number of permutations of length $n$ is $n!$, adding up,
and dividing by $n!$ gives a `combinatorial proof' that $c(n)=1$ for every $n$.
These coefficients also famously feature in
the {\it cycle index polynomial} of the symmetric group in Polya-Redfield theory. A much less famous occurrence,
kindly pointed to us by Andrew Sills, is due to Major Percy A. MacMahon ([M], p. 61ff), where these are
the coefficients in a certain `partial fraction decomposition', that has been beautifully extended by Sills([Si]) into
a multivariate identity.

Incidentally, the expression for $c(n)$ is precisely the {\it homotopy cardinality} $|X| = 1$ of the `action groupoid' $X = S_n\, //\, S_n$ whose objects are permutations of $n$ and whose isomorphisms are generated by the conjugation action of $S_n$ on itself.
Isomorphism classes of objects in $S_n\, //\, S_n$ are represented by partitions (since two permutations are conjugate just in case they have the same cycle lengths), and the formula $(1^{a_1} a_1!) \cdots (n^{a_n} a_n!)$ gives the order of the automorphism group 
($\cong$ stabilizer subgroup for the conjugation action)
of any permutation with cycle structure $1^{a_1} \dots n^{a_n}$.

{\bf Another  easy question}

If you look at partitions in {\bf frequency notation} $1^{a_1} \dots n^{a_n}$ and now give each of them `credit'
$$
\frac{1}{ (1!^{a_1} a_1!) \cdots (n!^{a_n} a_n!)} \quad ,
$$ 
and define $d(n)$ as the sum of these over all partitions of $n$, then
we have
$$
\sum_{n=0}^{\infty} d(n) q^n = e^{q/1!} \, e^{q^2/2!} \, e^{q^3/3!} \cdots = exp\, \left ( \sum_{i=1}^{\infty} \frac{q^i}{i!} \right )
\, = \,
exp(exp(q)-1))= \sum_{n=0}^{\infty} \frac{B_n}{n!} q^n \quad ,
$$
where $B_n$ are the Bell numbers, {\tt https://oeis.org/A000110}. So $d(n)=\frac{B_n}{n!}$. In particular,
$$
\sum_{n=0}^{\infty} d(n) \, =\, e^{e-1} \,= \,      5.5749415247608806\dots
$$
is the homotopy cardinality of the groupoid whose objects are finite sets equipped with a partition and whose isomorphisms are partition-respecting bijections.

{\bf The beautiful work of Robert Schneider}

While we believe that Questions 1 and 2 are new, it turns out that our kind of `fractional counting' of partitions 
showed up recently in
a remarkable PhD thesis, written by Robert Schneider ([Sc1], see also [Sc2]), where the sum of the reciprocals of
{\bf squares} (and more general powers)
of the `product of the parts', (that he calls the {\it norm}), taken over {\it all} partitions into even parts
(and more generally multiples of $m$, for any $m \geq 2$) are very elegantly expressed in terms of values of
the Riemann Zeta function at integer arguments.
We thank Andrew Sills for bringing this to our attention.

{\bf Epilogue: We need yet another On-Line Encyclopedia}

In addition to the great OEIS ([Sl]) that lets you identify integer sequences, and the very useful ``Inverse Symbolic
Calculator'' ([BP]) that lets you identify constants, it would be useful to have a searchable database of
continuous functions defined (for starters) on $0<x<1$, given numerically with, say, a resolution of $0.01$,
so each function will have $100$ values in floating point. If such a data-base existed, Question 2 may have been
answered (but one would need to go pretty far to get good approximations for $f(x)$).

\bigskip

\vfill\eject

{\bf Note Added Oct. 31, 2018}: A few hours after the posting of this article in the arxiv, Will Sawin, Columbia University,
solved Question 1. A donation to the OEIS, in his honor, has been made. Below is his argument.

\centerline
{\bf Appendix: Solution to Question 1 }
\bigskip
\centerline
{\it Will SAWIN}
\bigskip

In this appendix, we show that in the series

$$
\sum_{n=0}^{ \infty}\, b(n) q^n = \prod_{n=1}^{\infty   } \frac{1}{1-\frac{q^n}{n}} \quad .
$$

we have $$ \lim_{n \to \infty} \frac{b(n)}{n} = e^{-\gamma}$$

Let us choose $c(n)$ so that we have
$$
\sum_{n=0}^{ \infty} c(n) q^n = (1-q)^2 \prod_{n=1}^{\infty   } \frac{1}{1-\frac{q^n}{n}} \quad .
$$

We will show that $$ \sum_{n=0}^{\infty} |c(n)|$$ is bounded and that $$ \sum_{n=0}^{\infty} c(n) = e^{\gamma}.$$ The result will then follow because the coefficients of $\frac{1}{ 1-q} \sum_{n=0}^{\infty} c(n) q^n$ are the partial sums of $\sum_{n=0}^{\infty} c(n) q^n$ and so converge to $e^{-\gamma}$, and thus the coefficients of $\frac{1}{ (1-q)^2} \sum_{n=0}^{\infty} c(n) q^n$, which are the partial sums of those coefficients, are asymptotic to $e^{-\gamma} n$. 

To understand $c(n)$, using the formal identity $$
\frac{1}{1-q} = \prod_{n=1}^{\infty} e^{ \frac{q^n}{n} } \quad ,
$$ 
we have
$$
 (1-q)^2 \prod_{n=1}^{\infty   } \frac{1}{1-\frac{q^n}{n}}
=  (1-q) \prod_{n=2}^{\infty   } \frac{1}{1-\frac{q^n}{n}}
= e^{ -q} \prod_{n=2}^{ \infty} \frac{e^{ - \frac{q^n}{n}}}{1-\frac{q^n}{n}} \quad .
$$

Let us first evaluate this infinite product with $1$ substituted for $q$. This is
$$
 e^{ -1} \prod_{n=2}^{\infty} \frac{e^{ - 1/n}}{1-1 /n}
=  \lim_{m \rightarrow  \infty} \left [ \exp( \sum_{n=1}^m -\frac{1}{n} )  \prod_{n=2}^{m} \frac{1}{1-\frac{1}{n}} \right ]
 =  \lim_{m \rightarrow  \infty} \left [ \left ( \exp( \sum_{n=1}^m - \frac{1}{n} ) \right ) \cdot  m \right ]
$$
$$
 =   \lim_{m \rightarrow \infty} \exp(\log m -  \sum_{n=1}^m \frac{1}{n} ) 
  = e^{ - \gamma}  \quad .
$$

Now using the corollary 
$$
\frac{1}{1- \frac{q^n}{n}} = \prod_{d=1}^{\infty} e^{ (\frac{q^n}{n})^d/d }  \quad ,
$$
of the same formal identity, 
we see that the terms except $e^{-q}$ in the infinite product have nonnegative coefficients . Because we have already shown that the infinite product converges, it converges also if we take the absolute value of each coefficient, and so $ \sum_{n=0}^{\infty} |c(n)|$ is bounded (in fact by $e^{2 -\gamma}$ as taking absolute values of coefficients turns the $e^{-q}$ term into $e^q$). Because the product and sums converge absolutely, the partial sums converge to the limit of the product.

This research was conducted while Will Sawin was a Clay Research Fellow. \halmos

[Update 10:30am, Nov. 2, 2018: Fernando Chamizo has independently solved the first question, with an
even shorter rendition. It is available, with kind permission of Professor Chamizo, from \hfill\break
{\tt http://sites.math.rutgers.edu/\~{}zeilberg/mamarim/mamarimPDF/fcpFernandoChamizo.pdf } \quad .]

[Update 2:30pm, Nov. 2, 2018:  Gjergji Zaimi kindly pointed out that Question 1 was posed and answered
completely by D.H. Lehmer in Theorem 3 of his article  \hfill\break
( {\tt http://matwbn.icm.edu.pl/ksiazki/aa/aa21/aa21123.pdf} )
``On reciprocally weighted partitions", Acta Arithmetica XXI(1972), p. 379-388.]

[Update Nov. 6, 2018: Laurent Habsieger has independently solved the first question.
It is available, with kind permission of Professor Habsieger, from \hfill\break
{\tt http://sites.math.rutgers.edu/\~{}zeilberg/mamarim/mamarimPDF/fcpLaurentHabsieger.pdf } \quad .]

\bigskip

{\bf References}
\bigskip

[A]  George Andrews, ``{\it The Theory of Partitions}'', Cambridge University Press, 1984.
Originally published by Addison-Wesley, 1976.

[BD] John Baez and James Dolan, {\it From finite sets to Feynman diagrams}, in
B. Engquist and W. Schmid, editors,
{\it Mathematics Unlimited---2001 and Beyond}, 29--50, Springer, 2001.

[BP] Jonathan Borwein, Simon Plouffe et. al, {\it The Inverse Symbolic Calculator},  \hfill\break
{\tt https://isc.carma.newcastle.edu.au/} \quad  .

[FGZ]
Jane Friedman, Ira Gessel and Doron Zeilberger, {\it Talmudic lattice path counting},
J. Combinatorial Theory (ser. A) {\bf 68} (1994), 215--217. Available from \hfill\break
{\tt http://sites.math.rutgers.edu/\~{}zeilberg/mamarim/mamarimhtml/talmud.html} \quad .

[M] Percy Alexander MacMahon, {\it ``Combinatory Analysis, vol. II,''} Cambridge Univ. Press, 1916. 

[Sc1] Robert Schneider, {\it ``Eulerian series, zeta functions and the arithmetic of partitions''}, Ph.D. thesis, 2018. Available from \hfill\break
{\tt http://www.mathcs.emory.edu/\~{}rpschne/} [accessed Oct. 26, 2018] \quad .

[Sc2] Robert Schneider, {\it Partition Zeta Functions}, Research in Number Theory (2016) {\bf  2:9}.

[Si] Andrew Sills, {\it The combinatorics of MacMahon's partial fractions},  Annals of Combinatorics, to appear. Available from \hfill\break
{\tt http://home.dimacs.rutgers.edu/\~{}asills/MacMahon/MacMahon2017.pdf} \quad .

[Sl] Neil James Alexander Sloane, {\it The On-Line Encyclopedia of Integer Sequences}, \hfill\break
{\tt https://oeis.org} \quad .
\bigskip
\hrule
\bigskip
Doron Zeilberger, Department of Mathematics, Rutgers University (New Brunswick), Hill Center-Busch Campus, 110 Frelinghuysen
Rd., Piscataway, NJ 08854-8019, USA. \hfill\break
Email: {\tt DoronZeil@gmail.com}   \quad .
\bigskip
Noam Zeilberger, School of Computer Science, University of Birmingham, Birmingham, UK. \hfill\break
Email: {\tt  noam.zeilberger@gmail.com}   \quad .
\bigskip

First Written: Oct. 30, 2018. This version: Nov. 6, 2018.
\end